\documentclass[10pt,shortpaper,twoside]{IEEEtran}
\usepackage{cite}
\usepackage{amsthm}
\usepackage{amsmath,amsfonts}
\usepackage{algorithmic}
\usepackage{siunitx}
\usepackage{graphicx}
\usepackage{widetext}
\usepackage{balance}
\usepackage{hyperref}
\usepackage{xcolor}
\usepackage{wrapfig}
\usepackage[left=48pt, right=48pt, bottom=44pt, top=60pt]{geometry}
\def\BibTeX{{\rm B\kern-.05em{\sc i\kern-.025em b}\kern-.08em
		T\kern-.1667em\lower.7ex\hbox{E}\kern-.125emX}}
	
\usepackage{amssymb,verbatim,bm,color,times,mathtools,multicol,cuted}

\newcommand{\trace}{\mathrm{tr}}
\newcommand{\matrixstyle}[1]{\mathrm{#1}}

\DeclareMathOperator{\expect}{\mathbb{E}}

\newtheorem{assumption}{\textbf{Assumption}}
\newtheorem{definition}{\textbf{Definition}}

\begin{document}
	\title{\Huge Risk Sensitive Path Integral Control for \\ Infinite Horizon Problem Formulations}
	\author{
		Tom Lefebvre and Guillaume Crevecoeur
		\thanks{\{Tom.Lefebvre,Guillaume.Crevecoeur\}@ugent.be are with the Dept. of Electromechanical, Systems and Metal Engineering, Ghent University, 9000 Ghent; and with the Core Lab EEDT-DC, Flanders Make.
		}
	\thanks{Manuscript received Month XX, 20XX; revised Month XX, 20XX;
		accepted Month XX, 20XX. Date of publication Month XX, 20XX;
		date of current version  Month XX, 20XX. Recommended by Senior
		Editor YY}
	\thanks{Digital Object Identifier ZZZ}
		\vspace*{-12pt}
	}

	\maketitle

	\begin{abstract} 
	Path Integral Control methods were developed for stochastic optimal control covering a wide class of finite horizon formulations with control affine nonlinear dynamics. Characteristic for this class is that the HJB equation is linear and consequently the value function can be expressed as a conditional expectation of the exponentially weighted cost-to-go evaluated over trajectories with uncontrolled system dynamics, hence the name. Subsequently it was shown that under the same assumptions Path Integral Control generalises to finite horizon risk sensitive stochastic optimal control problems. 
	Here we study whether the HJB of infinite horizon formulations can be made linear as well. Our interest in infinite horizon formulations is motivated by the stationarity of the associated value function and their inherent dynamic stability seeking nature. Technically a stationary value function may ease the solution of the associated linear HJB. Second we argue this may offer an interesting starting point for off-linear Reinforcement Learning applications. We show formally that the discounted and average cost formulations are respectively intractable and tractable. 
	\end{abstract}
	

\section{Introduction}
\IEEEPARstart{R}{isk} Sensitive Optimal Control (RSOC) is a generalisation of Stochastic Optimal Control (SOC) that takes into account the risk seeking or risk averse nature of the controller \cite{fleming1992risk}. As opposed to the expected cost-to-go performance criteria in SOC, RSOC considers an expected exponentiated cost-to-go criteria. Depending on a risk sensitivity parameter $\phi$, the paradigm puts emphasis on the mode either the tail of the cost-to-go distribution which may be understood as emphasizing the risk-averse or -seeking incentive of the resulting optimal controller. The risk neutral setting where $\phi = 0$ recovers the conventional paradigm. RSOC was founded first in discrete time \cite{howard1972risk} and was generalised the following year to continuous time for linear systems \cite{jacobson1973optimal}. Both considered finite horizons. RSOC has been well studied since and has been generalised to discounted and average cost infinite horizons \cite{fleming1995risk,bensoussan1992bellman}. 

On account of Bellman's principle of Dynamic Programming, the value function and optimal control must satisfy a second order nonlinear Partial Differential Equation (PDE), the so called Hamilton-Jacobi-Bellman (HJB) equation (or Bellman recurrence equation for discrete time systems). The benefit of Dynamic Programming is the closed-loop or feedback control (i.e. the optimal control is a function of the state). Solving the HJB equation for arbitrary process dynamics and cost would thus produce the global optimal feedback control that is sought for in many applications. Despite the elegant theory, however, exact solutions of the HJB are rare and restricted to linear systems with quadratic cost formulations. If the problem formulation does not satisfy these conditions, one has to rely on approximation techniques \cite{tassa2007lshjb} or be satisfied with local solution that can be recovered using gradient based optimization \cite{mayne1966ddp} or trajectory optimization methods \cite{tassa2012synthesis}.

Recently it was shown that under certain assumptions the HJB becomes linear. Using the Feynman-Kac formula \cite{feynman2005space} the optimal value function can be expressed as the $\log$-transform of the expected exponentiated cost-to-go taken over trajectories with uncontrolled process dynamics. This amounts to a path integral formalism. The pioneering work was done by Kappen for finite horizon risk neutral formulations in \cite{kappen2005linear} and was generalised subsequently for finite horizon risk sensitive formulations in \cite{vdb2010risk}. The first application of the $log$-transform to stochastic optimal control is due to Fleming \cite{fleming1977exit} and Holland \cite{holland1977new}. However, Kappen was first to realize that control laws could be computed based on algorithms for probabilistic inference giving rise to Path Integral Control (PIC) methods. Numerous algorithms and applications have been documented since ranging from robot dog control \cite{theodorou2010reinforcement}, model based path integral control \cite{williams2017model} and reinforcement learning \cite{thalmeier2020adaptive}.

A similar path has been explored for Linearly Solvable Optimal Control (LSOC), sometimes referred to as KL-divergence control \cite{Dvijotham}. LSOC is a discrete time control framework that is closely related to PIC. Here the trick involves lifting the optimization entirely to the state space and considering a controlled state transition probability distribution \cite{lefebvre2020elsoc}. A cost-to-go is defined in terms of a state dependent cost rate and a weighted KL-divergence between the controlled and uncontrolled transition distribution. A unified framework has been proposed in \cite{Dvijotham} covering risk sensitive finite horizon and discounted and average cost infinite horizon formulations.

As far as we know of Risk Sensitive PIC (RSPIC) has not been generalised to infinite horizon formulations despite that for many applications an infinite horizon is natural. E.g. the discounted infinite horizon is the overarching framework of many Reinforcement Learning algorithms. Second we argue that the stationarity of the associated value functions may ease the solution of the associated linear HJB \cite{tassa2007lshjb}. Finally we reason the infinite horizon RSPIC formulations may offer an interesting starting point for off-linear RL applications where the main goal is to derive an optimal policy from historical data \cite{agarwal2020optimistic,levine2020offline}. By definition the historical data must be collected from the uncontrolled or suboptimally stabilized system linking directly with the expressions in RSPIC.

We show formally that discounted and average cost RSPIC formulations are respectively intractable and tractable. In particular, the $\phi$-sensitive value function is related to the $(\phi-\psi)$)-sensitive value function evaluated with uncontrolled dynamics. For discounted costs, the $(\phi-\psi)$)-sensitive value function gives a lower bound on the $\phi$-sensitive value function. For average costs, they are equal. 

\section{Preliminary results}
\begin{figure*}
	\begin{equation}
	\label{eq:1}
	\begin{aligned}
	Z(A,x) &= \expect_x\left[\exp\left(A\int_0^r e^{-\alpha t} q(X(t))\text{d}t\right)\exp\left(A e^{-\alpha r}\int_r^\infty e^{-\alpha t} q(X(t))\text{d}t\right)\right]  \\
	&=  \expect_x\left[\exp\left(A\int_0^r e^{-\alpha t} q(X(t))\text{d}t\right)Z\left(A e^{-\alpha r},X(r)\right)\right]
	\end{aligned}
	\end{equation}
	\begin{equation}
	\label{eq:2}
	\begin{aligned}
	Z(A e^{-\alpha r},x) &= \expect_x\left[\exp\left(A e^{-\alpha r}\int_0^\infty e^{-\alpha t} q(X(t))\text{d}t\right)\right] = \expect_x\left[\exp\left(A e^{-\alpha B r}\int_0^\infty e^{-\alpha t} q(X(t))\text{d}t\right)^{e^{-\alpha(1-B) r}}\right] \\
	&\leq Z(Ae^{-\alpha B r} ,x)^{e^{-\alpha(1-B) r}}
	\end{aligned}
	\end{equation}
	\begin{equation}
	\label{eq:3}
	Z(A,x) \leq \expect_x\left[\exp\left(A\int_0^r e^{-\alpha t} q(X(t))\text{d}t\right)Z\left(A e^{-\alpha B r} ,X(r)\right)^{e^{-\alpha (1-B)r}}\right]
	\end{equation}
	\hrulefill
\end{figure*}
Before we discuss RSOC and RSPIC we formally study the behaviour of certain infinite horizon functionals under the following stochastic process model. These infinite integrals will emerge when we discuss infinite horizon Path Integral Control formulations and in this context allows for fast interpretation of the result.

We consider a stochastic process model given by the following It\^{o} equation
\begin{align*}
\text{d}X &= a(X(t))\text{d}t + \sigma(X(t)) \text{d}W(t) \\
X(0) &= x
\end{align*}
Here $x\in\mathbb{R}^n$ is the state, $W$ is a Wiener process in $\mathbb{R}^p$ which models the noise in the system, $a$ represents the time invariant autonomous dynamics and $\sigma$ is a constant full rank matrix with appropriate dimensions.

\subsection{Discounted cost}
First consider the integral with constant $A \neq 0$, $\alpha > 0$ and where $q$ is a function defined on $\mathbb{R}\times\mathbb{R}^n$
\begin{equation*}
Z(A,x) = \expect_x\left[\exp\left(A\int_0^\infty e^{-\alpha t} q(X(t))\text{d}t\right)\right]
\end{equation*}

This expression represents a solution of a stochastic Partial Differential Equation (PDE). To see this we follow the argument made in equation (\ref{eq:1}). Using It\^{o}'s lemma for $Z(A e^{-\alpha r},X(r))$ and $r\rightarrow 0$ we obtain
\begin{equation*}
\alpha A \partial_\phi Z = A q Z + a \cdot \nabla Z + \tfrac{1}{2} \Delta_\Sigma Z 
\end{equation*}
where $\Sigma = \sigma\sigma^\top$ and $\Delta_\Sigma Z = \trace (\Sigma \nabla Z)$.

This formal derivation suggests that $Z(A,x)$ is an exact solution of the differential equation above. The existence and uniqueness are not treated but we reason they can be verified using similar methods of proof as in \cite{menaldi2005remarks} Theorem 2.1.

Second, Jensen's inequality implies the inequality in equation (\ref{eq:2}) with $B < 1$. Substitution of this expression for $Z(A e^{-\alpha r},X(r))$ in equation (\ref{eq:1}) generates the inequality presented in equation (\ref{eq:3}). Again using It\^{o}'s lemma for $Z(Ae^{-\alpha B r} ,X(r))^{e^{-\alpha (1-B) r}}$ and $r\rightarrow 0$ we can obtain the following Partial Differential inequality
\begin{equation}
\label{eq:if}
\begin{aligned}
\alpha A B \partial_A Z + \alpha (1-B) Z \log Z &\leq A q Z + a \cdot \nabla Z + \tfrac{1}{2} \Delta_\Sigma Z 
\end{aligned}
\end{equation}
For $B> 1$ the inequality flips and for $B=1$ clearly the equality holds.

Put differently, when $B>1$ then $Z(A,x)$ gives a lower bound for the solution of the associated PDE otherwise for $B< 1$ expression $Z(A,x)$ gives an upper bound. For $B= 1$ we recover the PDE in (\ref{eq:if}).

\subsection{Average cost}

Second consider the limit with $A \neq 0$
\begin{equation*}
\chi(A) = \lim_{T\rightarrow \infty} \tfrac{1}{T} \log \expect_{x} \left[\exp\left(\int_0^T q(X(t))\text{d}t\right)^A\right]^{\tfrac{1}{A}} 
\end{equation*}
Given assumptions for which the previous limit exists (see section \ref{sec:average-cost-rsoc}), next we consider the integral expression
\begin{equation*}
z(A,x) = \expect_{x}\left[\exp\left(A\int_0^\infty (q(X(t))-\chi(A))\text{d}t\right)\right]
\end{equation*}
Given assumptions for which this expression converges (again see section \ref{sec:average-cost-rsoc}), $z$ represents a solution to a stochastic PDE. To see this we argue as follows
\begin{align*}
z(A,x) = \expect_{x}\left[\exp\left(A\int_0^r (q(X(t))-\chi(A))\text{d}t\right)z(A,X(r))\right]
\end{align*}
Again using It\^{o}'s lemma for $Z(A,X(r))$ and $r\rightarrow 0$ we obtain the following PDE
\begin{equation*}
0 = A (q-\chi(A)) Z + a\cdot \nabla z + \tfrac{1}{2}\Delta_\Sigma z
\end{equation*}

Again this suggests that $z(A,x)$ is a solution of the PDE. The existence of $z(A,x)$ will follow from our subsequent discussion.

\section{Risk Sensitive Optimal Control}
In this section we move to RSOC in continuous time. Therefore let us first extend our dynamic system model to controlled dynamics
\begin{align*}
\text{d}X(t) &= f(t,X(t),u(t))\text{d}t + \sigma(t,X(t))\text{d}W(t) \\
X(0) &= x
\end{align*}
where $u\in\mathbb{R}^m$ is the control.

Next we study some properties of risk sensitive optimal control formulations. We treat the finite horizon setting as a point of reference. Then we introduce the discounted infinite horizon and average cost infinite horizon formulations.

\subsection{Finite Horizon RSOC} 
In a finite horizon setting we evaluate the performance of the system using the cost
\begin{equation*}
J_T^u(t,x)[l,m] = \int_{t}^{T} l(t,x(t),u(t))\text{d}t + m(x(T))
\end{equation*}
where $l$ is an instantaneous cost and $m$ is the final cost. For some parameter $\phi$ we define the risk sensitive cost-to-go
\begin{equation*}
V^u_T(\phi,t,x) = \left\lbrace \begin{aligned}
&\log \expect^u_{t,x}[\exp\left( J_T^u(t,x)\right)^\phi]^{\tfrac{1}{\phi}}, && \phi \neq 0 \\
&\expect^u_{t,x}[J_T^u(t,x)], && \phi = 0
\end{aligned} \right.
\end{equation*}
Here $\expect^u_{t,x}$ denotes expectation over all realisations of
the dynamics starting in $x$ at time $t$ with control $u$.

\begin{definition}
\label{def:rsfh}
The finite horizon risk sensitive value function is defined as
\begin{equation*}
V_T(\phi,t,x) = \min_u V^u_T (\phi,t,x)
\end{equation*}
\end{definition}

The control $u^*$ satisfying $V_T = V_T^{u^*}$ is called optimal. In this setting $\phi$ expresses the risk sensitivity. When $\phi > 0$, the controller is said to be risk averse, when $\phi < 0$, the controller is said to be risk seeking, in the limit $\phi \rightarrow 0$ the problem reduces to the risk neutral or standard formulation. To see this note that when $\phi$ is small we have
\begin{equation*}
V_T^u= \expect[J^u_T] + \tfrac{1}{2} \phi\left( \expect[(J^{u}_T)^2] -  \expect[J^{u}_T]^2 \right)+ \mathcal{O}(\phi^2)
\end{equation*}

The function $V_T$ satisfies the following HJB equation.  

\begin{definition} The finite horizon risk sensitive HJB equation is defined as 
\begin{align*}
-\partial_t V_T &=  \min_u \left\{l + f \cdot \nabla V\right\} + \tfrac{1}{2}\Delta_\Sigma V_T + \tfrac{1}{2}\phi\|\nabla V_T\|^2_\Sigma \\
V_T(x) &= m(x)
\end{align*}
\end{definition} 

We give a formal derivation. We start with 
\begin{equation*}
W_T = \min_u \exp(\phi V^u_T)
\end{equation*}
Then for any $r>t$ and for any control $u(t) = u(t,x(t))$ we give the argument in equation (\ref{eq:4}). Using It\^{o}'s lemma for $W_T(\phi,r,X(r))$ and $r\rightarrow t$ we obtain
\begin{equation*}
\partial_t W_T  + \tfrac{1}{2} \Delta_\Sigma W_T + \min_u \left\{\phi l W_T + f \cdot \nabla W_T\right\} = 0
\end{equation*}
Next for $\phi > 0$ we set $W_T = \exp(\phi V_T)$ and deduce that 
\begin{equation*}
-\partial_t V_T =  \min_u \left\{l + f \cdot \nabla V_T\right\} + \tfrac{1}{2}\Delta_\Sigma V_T + \tfrac{1}{2}\phi\|\nabla V_T\|^2_\Sigma 
\end{equation*}
For $\phi < 0 $, $W_T$ is defined as the maximum, otherwise the derivation is similar. For $\phi = 0$ we recover the standard finite horizon HJB.

In the general case this HJB is intractable.


\subsection{Discounted Infinite Horizon RSOC} 

Now we treat the discounted infinite horizon setting. Here we evaluate the performance of the system using the cost
\begin{equation*}
J_\alpha^u(x)[l] = \int_{0}^{\infty} e^{-\alpha t}l(x(t),u(t))\text{d}t
\end{equation*}
where $l$ is a time invariant cost and $\alpha > 0$ is the discount factor. In this setting we further assume time invariant dynamics. For some parameter $\phi$ we define the risk sensitive cost-to go
\begin{equation*}
V_\alpha^u (\phi,x) =\log\expect^u_{x}[\exp\left(J_\alpha^u(x)\right)^\phi]^{\frac{1}{\phi}}
\end{equation*}
\begin{definition}
	\label{def:rsdih}
	The discounted infinite horizon risk sensitive value function is defined as
	\begin{equation*}
	V_\alpha(\phi,x) = \min_u V_\alpha^u (\phi,x)
	\end{equation*}
\end{definition}

The stationary policy $u^*$ satisfying $V_\alpha = V_\alpha^{u^*}$ is called optimal. The function $V_\alpha$ is also defined by a HJB. 
\begin{definition} The discounted infinite horizon risk sensitive HJB equation is defined as 
	\begin{equation*}
	\alpha\left(V_\alpha + \phi \partial_\phi V_\alpha \right) = \min_u \left\{l + f \cdot \nabla V_\alpha\right\} + \tfrac{1}{2}\Delta_\Sigma V_\alpha + \tfrac{1}{2}\phi\|\nabla V_\alpha\|^2_\Sigma 
	\end{equation*}
\end{definition} 

Formally, for $\phi > 0 $ we define 
\begin{equation*}
W_\alpha = \min_u \exp(\phi V_\alpha^u)
\end{equation*}
Then for any $r>0$ and for any control $u(t) = u(x(t))$ we follow the argument in equation (\ref{eq:5}). Using It\^{o}'s lemma for $W_\alpha(\phi e^{-\alpha r},X(r))$ and $r\rightarrow 0$ we obtain
\begin{equation*}
-\alpha \phi \partial_\phi W_\alpha  + \tfrac{1}{2} \Delta_\Sigma W_\alpha + \min_u \left\{\phi l W_\alpha + f \cdot \nabla W\right\} = 0
\end{equation*}
Next we set $W_\alpha = \exp(\phi V_\alpha)$ and deduce that 
\begin{equation*}
\alpha\left(V_\alpha + \phi \partial_\phi V_\alpha \right) = \min_u \left\{l + f \cdot \nabla V_\alpha\right\} + \tfrac{1}{2}\Delta_\Sigma V_\alpha + \tfrac{1}{2}\phi\|\nabla V_\alpha\|^2_\Sigma 
\end{equation*}
For $\phi < 0 $, $W_\alpha$ is defined as the maximum, otherwise the derivation is similar and amounts to the same result. For $\phi = 0$ we recover the standard discounted infinite horizon HJB.

\begin{figure*}
	\begin{equation}
	\label{eq:4}
	\begin{aligned}
	W_T(\phi,t,x) &= \min_u \expect_{t,x} \left[\exp\left(\phi \int_t^T l(\tau,X(\tau),u(\tau))\text{d}\tau + \phi m(X(T))\right)\right] \\
	&= \min_u \expect_{t,x} \left[\exp\left(\phi \int_t^r l(\tau,X(\tau),u(\tau))\text{d}\tau \right) W_T(\phi,r,X(r))\right]
	\end{aligned}
	\end{equation}
	\begin{equation}
	\label{eq:5}
	\begin{aligned}
	W_\alpha(\phi,x) &= \min_u \expect_{x} \left[\exp\left(\phi \int_0^\infty e^{-\alpha \tau}l(X(\tau),u(\tau))\text{d}\tau \right)\right] \\
	&= \min_u \expect_{x} \left[\exp\left(\phi \int_0^r e^{-\alpha \tau} l(X(\tau),u(\tau))\text{d}\tau \right) W_\alpha(\phi e^{-\alpha r},X(r))\right]
	\end{aligned}
	\end{equation}
	\hrulefill
\end{figure*}
\subsection{Average Cost RSOC}\label{sec:average-cost-rsoc}
Finally, let us treat the average cost setting. The formal treatment is less straightforward then the finite and discounted counterparts. It can be derived as a limit case of the finite horizon setting ($T\rightarrow \infty$) \cite{nagai1995bellman} or as a limit case of the infinite horizon setting ($\alpha \rightarrow 0$) \cite{menaldi2005remarks}. Here we opt for the finite treatment. Following the argument in the finite horizon setting we have that
\begin{equation*}
-\partial_t V_T = \min_u \left\{l+ f\cdot \nabla V_T \right\} + \tfrac{1}{2}\Delta_\Sigma V_T + \tfrac{1}{2} \phi \|\nabla V_T\|^2_\Sigma
\end{equation*}

Given that $f$ and $l$ are periodic in $x$, $f$ is continuous in $x$ and $u$ and Lipschitz continuous in $x$ and additionally $l$ is continuous and non-negative, it can be shown that the following limit exists \cite{nagai1995bellman} (or \cite{menaldi2005remarks} for a similar result but starting from the discounted setting).
\begin{equation*}
\chi = \lim_{T\rightarrow \infty} \tfrac{1}{T} V_T(\phi,0,x) = \lim_{T\rightarrow \infty} -\partial_t V_T(\phi,0,x)
\end{equation*}
Here $\chi$ is referred to as the average cost per stage and can be understood as the minimal asymptotic cost accumulation rate at large time horizons.

In this context we can not define a value function in the same way as we did in the finite and discounted settings. In the average cost setting an value function $v$ is defined that has the meaning of a differential value function. The choice of the reference state is arbitrary and will only change the value of $v$ with some constant. Here we take the origin.
\begin{definition} \label{def:rsacih}
	The average cost infinite horizon risk sensitive value function is defined as
	\begin{align*}
	v(\phi,x) &= \lim_{T\rightarrow \infty} V_T(\phi,0,x) - V_T(\phi,0,0) 
	\end{align*}
\end{definition}
Although this result is not intuitive, one might understand this as follows. Assume that we follow two particles, initialised at $x$ and $0$ and measure there accumulated cost-to-go under the optimal policy. Although each of our cost measurements will converge to $\chi T$ as time passes, if we would follow these particles long enough, we would notice that the difference in cost measurements becomes constant and depends only on the relative difference between the states the particles were initialized hence the existence of a non-constant $v(\phi,x)$.

The following HJB relates $\chi$ and $v$.
\begin{definition} The average cost infinite horizon risk sensitive HJB equation is defined as 
	\begin{equation*}
	\chi = \min_u \left\{l + f \cdot \nabla v \right\} + \tfrac{1}{2}\Delta_\Sigma v + \tfrac{1}{2}\phi\|\nabla v\|^2_\Sigma 
	\end{equation*}
\end{definition}

Note that for the optimal policy $u^*(\phi,x)$ the HJB implies
\begin{equation*}
v(\phi,x) = \log \expect_x\left[\exp\left(J_\infty^{u^*}(0,x)[l-\chi,0]\right)^{\phi}\right]^{\tfrac{1}{\phi}}
\end{equation*}
The formal derivation is similar as for the average cost in the preliminary section substituting $v = \phi \log z$. The intuitive understanding is that if we subtract the average cost per stage $\chi$ from the cost rate $l$, on large time horizons the difference between $\chi T$ and the accumulated cost is determined by the state the system was initialised. This again suggests that the average cost formulation seeks for controllers that initiate a stable limit cycles as are sought for in robotic locomotion as well as industrial applications.

\section{Risk Sensitive Path Integral Control}
In this section we show that for certain assumptions the HJB becomes linear and can be solved using the Feynman-Kac formula in the finite horizon setting. We treat this setting first. Then we move to discounted infinite horizon and illustrate why in this setting a linear solution does not exist. In conclusion we treat average cost problem formulations. As was mentioned in the introduction, we are motivated to study infinite horizon formulation because here the value functions are stationary. We argue that this condition may ease the solution of the linear HJB, if it exists, and might give rise to new learning algorithms.

\subsection{Finite Horizon RSPIC} By introducing a number of assumptions on the generic finite horizon problem formulation we can arrive at a linear HJB equation which solution is given by the Feynman-Kac formula. This is the subject of the present section.

\begin{assumption}\label{as:1} The instantaneous cost rate is control quadratic $
l(t,x,u) = q(t,x) + \tfrac{1}{2}\|u\|^2_{\matrixstyle{R}(t)}$ where $\matrixstyle{R} \succ 0$.
\end{assumption}

\begin{assumption}\label{as:2}  The dynamics are control affine $f(t,x,u) = a(t,x) + \matrixstyle{B}(t,x) u$.
\end{assumption}

Under these assumptions the HJB gives the following optimal control
\begin{equation*}
u^*(\phi,t,x) = - \matrixstyle{R}(t)^{-1} \matrixstyle{B}(t,x) \nabla V_T(\phi,t,x)
\end{equation*}	
Substitution of the optimal control in the HJB yields
\begin{equation*}
-\partial_t V_T = q + a \cdot \nabla V_T + \tfrac{1}{2}\Delta_\Sigma V_T - \tfrac{1}{2}\|\nabla V_T\|^2_{\matrixstyle{A}}  + \tfrac{1}{2}\phi\|\nabla V_T\|^2_\Sigma 
\end{equation*}
where $\matrixstyle{A} = \matrixstyle{B}\matrixstyle{R}^{-1}\matrixstyle{B}^\top$. Next we set $V_T = -\lambda \log Z_T$ to find	
\begin{equation*}
\partial_t Z_T = -\tfrac{1}{\lambda}q Z_T + a \cdot \nabla Z_T + \tfrac{1}{2}\Delta_\Sigma Z_T + \tfrac{1}{2}\tfrac{1}{Z_T}\|\nabla Z_T\|^2_{\matrixstyle{M}} 
\end{equation*}
where $\matrixstyle{M} = \lambda \matrixstyle{A} - (1 + \lambda \phi ) \Sigma$. 

The final term prevents us from applying the Feynman-Kac formula. Thus we introduce a third assumption.
\begin{assumption}\label{as:3}  The noise and control cost rate are proportional $\matrixstyle{A} = \psi \Sigma$ with $\psi>0$.
\end{assumption}
Then if we choose $\lambda = \frac{1}{\psi -\phi}$ the linear HJB follows
\begin{definition} The finite horizon risk sensitive linear HJB is defined as 
\begin{equation*}
\partial_t Z_T = (\phi-\psi)q Z_T + a \cdot \nabla Z_T + \tfrac{1}{2}\Delta_\Sigma Z_T
\end{equation*}
\end{definition}
Then according to the Feynman-Kac formula we have
\begin{equation*}
Z_T(\phi,t,x) = \expect^0_{t,x}\left[\exp(J_T^0(t,x))^{\phi-\psi}\right]
\end{equation*}
Note that the expectation is taken over \textit{uncontrolled} trajectories. Finally substituting $Z_T$ as a function of $V_T$
\begin{align*}
V_T(\phi,t,x) &= \min_u V_T^u(\phi,t,x) \\
&= \log \expect^0_{t,x}\left[\exp(J^0_T(t,x))^{\phi-\psi}\right]^{\tfrac{1}{\phi - \psi}} \\
&= V^0_T(\phi-\psi,t,x)
\end{align*}

Put differently, under assumptions \ref{as:1}, \ref{as:2} and \ref{as:3}, the $\phi$-sensitive \textit{optimal} value function is equal to the $(\phi-\psi)$-sensitive uncontrolled cost-to-go thus implying that the $\phi$-sensitive optimal control is as risk seeking as the uncontrolled dynamics in the $(\phi-\psi)$-sensitive setting. This result was presented earlier in \cite{vdb2010risk}. The relation between the $\phi$-sensitive optimal cost and $(\phi-\psi)$-sensitive uncontrolled cost-to-go and hence the $\phi$-sensitive optimal control and $(\phi-\psi)$-sensitive uncontrolled dynamics was not emphasized. In conclusion we emphasize that this does not imply that the $\phi$-sensitive optimally controlled dynamics and the $(\phi-\psi)$-sensitive uncontrolled dynamics are equivalent, it is only so that their averaged behaviour is similar.

\subsection{Discounted Infinite Horizon RSPIC} 

Here it is attempted to obtain a linear solution to the discounted infinite horizon HJB. 

Under assumptions \ref{as:1}, \ref{as:2} and \ref{as:3}, the optimal control is given
\begin{equation*}
u^*(\phi,x) = - \matrixstyle{R}^{-1} \matrixstyle{B}(x) \nabla V_\alpha(\phi,x)
\end{equation*}
Substitution into the HJB gives
\begin{multline*}
\alpha\left(V_\alpha + \phi \partial_\phi V_\alpha \right) = \\ q + a \cdot \nabla V_\alpha + \tfrac{1}{2}\Delta_\Sigma V_\alpha - \tfrac{1}{2}\|\nabla V_\alpha\|^2_{\matrixstyle{A}}  + \tfrac{1}{2}\phi\|\nabla V_\alpha\|^2_\Sigma 
\end{multline*}
Next we set
\begin{equation*}
V_\alpha = -\lambda \log Z_\alpha
\end{equation*}
Note that here $\lambda$ might depend on $\phi$ so that
\begin{align*}
\partial_\phi V_\alpha &= -\partial_\phi \lambda \log Z_\alpha - \tfrac{1}{Z_\alpha}\lambda \partial_\phi Z_\alpha \\
\nabla V_\alpha &= -  \tfrac{1}{Z_\alpha}\lambda\nabla Z_\alpha \\
\Delta_\Sigma V_\alpha &=  -  \tfrac{1}{Z_\alpha}\Delta_\Sigma Z_\alpha + \tfrac{1}{Z_\alpha} \lambda \|\nabla Z_\alpha\|^2_\Sigma
\end{align*}
We find 
\begin{multline*}
\alpha\phi \partial_\phi Z_\alpha + \alpha\left(1 + \tfrac{\phi}{\lambda}\partial_\phi \lambda \right) Z_\alpha \log Z_\alpha = \\
-\tfrac{1}{\lambda}q Z_\alpha + a \cdot \nabla Z_\alpha + \tfrac{1}{2}\Delta_\Sigma Z_\alpha + \tfrac{1}{2}\tfrac{1}{Z_\alpha}\|\nabla Z_\alpha\|^2_{\matrixstyle{M}} 
\end{multline*}
To get rid of the final term again we choose $\lambda = \frac{1}{\psi-\phi}$ so that $\partial_\phi \lambda = \frac{1}{(\psi-\phi)^2}$ which give the following linear HJB.
\begin{definition} The discounted infinite horizon risk sensitive linear HJB is defined as 
\begin{multline*}
\alpha\phi \partial_\phi Z_\alpha + \alpha\left(1 + \tfrac{\phi}{\psi-\phi} \right) Z_\alpha \log Z_\alpha = \\
(\phi-\psi) q Z_\alpha + a \cdot \nabla Z_\alpha + \tfrac{1}{2}\Delta_\Sigma Z_\alpha 
\end{multline*}
\end{definition}

For this framework to be useful a solution to the HJB must exist. However setting $A = \phi-\psi$ and $B = \frac{\phi}{\phi-\psi}$ and recalling the result from the introduction we find that for $B > 1$
\begin{equation*}
Z_\alpha(\phi,x) \geq Z_\alpha^0(\phi-\psi,x)
\end{equation*}
so that 
\begin{align*}
V_\alpha(\phi,x) &= \min_u V_\alpha^u(\phi,x) \\
&\geq \log \expect^0_{t,x}\left[\exp(J_\alpha^0(t,x))^{\phi-\psi}\right]^{\tfrac{1}{\phi - \psi}} \\
&= V_\alpha^0(\phi-\psi,t,x)
\end{align*}

It follows that it is impossible to establish a solution to the linear HJB in the discounted infinite horizon case. Hence neither we can estimate the discounted risk sensitive value function exactly. We can only obtain a lower bound. Since $V_\alpha < V_\alpha^u$ for any $u\neq u^*$, it follows that a lower bound on $V_\alpha$ is insufficient to derive a $\phi$-sensitive suboptimal control. We conclude that Path Integral Control does not generalise to discounted risk sensitive stochastic optimal control formulations.

\subsection{Average Cost RSPIC}

Finally, let us treat the average cost setting.
Under assumptions \ref{as:1}, \ref{as:2} and \ref{as:3}, again the value function is given by
\begin{equation*}
u^*(\phi,x) = - \matrixstyle{R}^{-1} \matrixstyle{B}(x) \nabla v(\phi,x)
\end{equation*}
Substitution in the HJB gives
\begin{equation*}
\chi = q + a \cdot \nabla v + \tfrac{1}{2}\Delta_\Sigma v - \tfrac{1}{2}\|\nabla v\|^2_{\matrixstyle{A}}  + \tfrac{1}{2}\phi\|\nabla v\|_\Sigma^2
\end{equation*}
Next we set $v = -\lambda \log z$ and find 
\begin{equation*}
0 =
-\tfrac{1}{\lambda}(q-\chi) z + a \cdot \nabla z + \tfrac{1}{2}\Delta_\Sigma z + \tfrac{1}{2}\tfrac{1}{z}\|\nabla z\|^2_{\matrixstyle{M}} 
\end{equation*}
Again we choose $\lambda = \frac{1}{\psi-\phi}$ to remove the final term and obtain the following linear HJB
\begin{definition} The average cost infinite horizon risk sensitive linear HJB is defined as 
\begin{equation*}
0 =
(\phi-\psi) (q - \chi) z + a \cdot \nabla z + \tfrac{1}{2} \Delta_\Sigma z
\end{equation*}
\end{definition} 

A solution to this HJB exists. Setting $A=\phi-\psi$ and recalling the result from the introduction we anticipate
\begin{equation*}
z(\phi,x) = \expect_x\left[\exp\left((\phi-\psi)\int_0^\infty  (q(X(t)-\chi))\text{d}t\right)\right]
\end{equation*}
where $\chi = \chi(\phi-\psi)$ and
\begin{equation*}
\chi(\phi-\psi) = \lim_{T\rightarrow \infty} \tfrac{1}{T} \log \expect_{x} \left[\exp\left(\int_0^T q(X(t))\text{d}t\right)^{\phi-\psi}\right]^{\tfrac{1}{\phi-\psi}} 
\end{equation*}

Substituting the expression as a function of $v$ for $z$ finally gives that
\begin{equation*}
v(\phi,x) =\log \expect_x\left[\exp\left(\int_0^\infty  (q(X(t)-\chi))\text{d}t\right)^{\phi-\psi}\right]^{\tfrac{1}{\phi-\psi}}
\end{equation*}

These results are confirmed by substituting the finite horizon Path Integral Control solution into the general expressions that generate the average cost setting by taking the large time limit of the finite horizon setting. In particular
\begin{align*}
\chi = \chi(\psi-\phi) &= \lim_{T\rightarrow \infty} \tfrac{1}{T} V_T(\phi,0,x) \\
&= \lim_{T\rightarrow \infty} \tfrac{1}{T} V_T^0(\phi-\psi,0,x)
\end{align*}
and
\begin{align*}
v(\phi,x) &= \lim_{T\rightarrow \infty} V^0_T(\phi-\psi,0,x) - V^0_T(\phi-\psi,0,0) 
\end{align*}

For the uniqueness and existence of a solution to the linear HJB we rely on the uniqueness and existence of a solution to the general nonlinear HJB.

\section{Example} As a running example we consider the linear system, $\text{d}X = (\matrixstyle{A} x + \matrixstyle{B} u)\text{d}t + \sigma \text{d}W$, quadratic cost rate $l=  \frac{1}{2}x^\top \matrixstyle{Q}x + \frac{1}{2} u^\top \matrixstyle{R} u$ then one may verify that 
\begin{align*}
z &= \exp\left((\phi-\psi)\tfrac{1}{2}x^\top \matrixstyle{S} x\right) \\
\chi &= \tfrac{1}{2}\trace (\matrixstyle{S}\Sigma)
\end{align*}
where $\matrixstyle{S}$ is the solution of the Ricatti equation
\begin{equation*}
0 = \matrixstyle{Q} + \matrixstyle{A}^\top \matrixstyle{S} + \matrixstyle{S}\matrixstyle{A} + (\phi-\psi)\matrixstyle{S}\Sigma \matrixstyle{S}
\end{equation*}
and if we substitute $\psi\Sigma = \matrixstyle{B}^\top \matrixstyle{R}^{-1} \matrixstyle{B}$ we recover 
\begin{equation*}
0 = \matrixstyle{Q} + \matrixstyle{A}^\top \matrixstyle{S} + \matrixstyle{S}\matrixstyle{A} + \matrixstyle{S}\left(\phi\Sigma - \matrixstyle{B}^\top \matrixstyle{R}^{-1}\matrixstyle{B}\right) \matrixstyle{S}
\end{equation*}
with $\phi = 0$ for example we retrieve the standard Ricatti equation of the continuous time Linear Quadratic Gaussian. For $\phi \neq 0$ we obtain the solution presented by Jacobson \cite{jacobson1973optimal}.

\section{Discussion}
The formal derivations presented in the previous section have generalised the Path Integral Control method to the average cost infinite horizon setting and illustrated why a generalisation to discounted infinite horizon is strictly impossible (under assumptions \ref{as:1}, \ref{as:2} and \ref{as:3}). However one has argued that average costs are in fact more realistic than discounted costs. In general, the average cost implies that some kind of dynamic stability of the system is looked for \cite{robin1983long}. Hence an average cost might be a natural choice for control tasks where stable limit cycles are sought such as industrial process control, robot locomotion, etc. It is rarely applied though since the associated HJB is difficult to solve specifically for large state-spaces. In their respect, the generalisation of PIC opens up interesting perspectives. Technically, the stationarity of the value function may ease direct solution of the value function.

An another important implication of our argument is that
\begin{equation*}
z(\phi,x) = C \expect_{x}\left[\exp\left(\int_0^Tq(X(t)\text{d}t\right)^{\psi-\phi} z(\phi,X(T))\right] 
\end{equation*}
where 
\begin{equation*}
C = e^{-\chi(\phi-\psi)T}
\end{equation*}
We argue that this stationary recurrence relation might be the starting point of off-line Reinforcement Learning algorithms where the main goal is to derive an optimal policy from historical data \cite{agarwal2020optimistic,levine2020offline}. Since here the expectation is taken inherently over uncontrolled system dynamics and the value function is stationary, a parametrised estimate might be obtained by approximating the expectation by means of Monte Carlo with the historical data and solving the linear equation. Clearly there are still several technical issues to overcome. For example the expectation is over trajectories and should be manipulated into an expectation over the stationary data set. Such an expression could then be used as the objective in a regression problem. Here the the risk sensitivity $\phi$ is an open tuning parameter. This is especially interesting for that given the same historical data set, a class of controllers could be obtained of varying robustness simply by varying $\phi$.

In conclusion, we highlight our main observations. We investigated Risk Sensitive Path Integral Control methods for infinite horizon problem formulations formally. The main goal is to verify under which assumptions a linear HJB equation can be derived for infinite horizon formulations. Here we were motivated to obtain a recurrence relation for the stationary value function that can be expressed as a conditional expectation over uncontrolled dynamics. Our formal derivations indicate that for discounted formulations this is impossible but for average cost formulations this is possible. We argue that the latter might offer an interesting starting point to derive offline Reinforcement Learning algorithms that can estimate a parametrised value function from historical data.

\section*{Acknowledgements} This research received funding from the Flemish Government under the ``Onderzoeksprogramma Artifici\"{e}le Intelligentie (AI) Vlaanderen'' programme.

\bibliographystyle{unsrt}
\bibliography{references}
	
\end{document}